# OPTIMIZING LOW-THRUST AND GRAVITY ASSIST MANEUVERS TO DESIGN INTERPLANETARY TRAJECTORIES[#]

**Massimiliano Vasile and Franco Bernelli-Zazzera**

*Dipartimento di Ingegneria Aerospaziale, Politecnico di Milano*
*Via La Masa 34, 20158, Milano, Italy*
*vasile@aero.polimi.it*

### Abstract

In this paper a direct method based on a transcription by Finite Elements in Time has been used to design optimal interplanetary trajectories, exploiting a combination of gravity assist maneuvers and low-thrust propulsion. A multiphase parametric approach has been used to introduce swing-bys among thrust and coast arcs. Gravity maneuvers are at first modeled with a link-conic approximation and then introduced through a full three-dimensional propagation including perturbations by the Sun. The method is successfully applied to the design of a mission to planet Mercury, for which different options corresponding to different sequences of gravity maneuvers or launch opportunities are presented.

## INTRODUCTION

Interplanetary missions pose extremely challenging design problems. In fact, with the exception of the planets closest to the Earth, a direct transfer requires a considerable amount of propulsion, which is best met by using high specific impulse propulsion systems like ion or plasma drives. Recently, the European Space Agency has analyzed the possibility to reach Mercury using Solar Electric Propulsion (SEP) as an alternative to chemical propulsion or multiple swing-bys[1]. However, for a mission to Mercury, the planet nearest to the Sun, the Δv requirements are so demanding that even SEP does not lead to a mission feasible in terms of mass consumption, duration and overall cost. The main problem

---

is the long thrust time required from a state-of-the-art low-thrust engine (over 7000 h of operations). An interesting way to reduce such a requirement is to resort to Gravity Assist (GA) maneuvers as shown by Langevin[2] in 1999. This concept requires nevertheless the need to establish the correct sequence of propulsion and GA maneuvers.

The design of a transfer trajectory combining SEP and GA can be regarded as a general trajectory optimization problem[3]. The dynamics of the spacecraft is governed mainly by the gravity attraction of the Sun, when the spacecraft is outside the sphere of influence of a planet, and by the gravity attraction of the planet during a gravity assist maneuver. Low-thrust propulsion is then used to shape trajectory arcs between two subsequent encounters and to meet the best incoming conditions for a swing-by.

The solution approach proposed by Langevin requires some heuristics and a remarkable amount of experience and physical insight into the problem. In fact thrust arcs and coast arcs are appropriately combined to shape correctly the trajectory. Most of the experience, however, is employed in designing appropriately the swing-bys guessing the right date for encounters and the correct set of parameters. Although this can be regarded as a direct optimization approach, leading to solutions really close to the optimum, it requires a considerable human effort.

A different option is to solve the problem through optimal control theory and Pontryagin maximum principle. As for all indirect approaches the solution found is extremely accurate but again it requires a considerable human effort in guessing correctly the switching structure and a first solution which leads the method to converge quickly to the optimum[4]. Examples of this approach can be found in the work of Colasurdo et al. or of Jehn and Kazkowsky[5].

An interesting alternative approach is to resort to direct collocation as demonstrated by Betts[6], who efficiently optimized a transfer trajectory to Mars combining low-thrust with two swing-bys of Venus.

In this paper an original direct optimization approach has been used to design an optimal interplanetary trajectory. The proposed approach is characterized by a transcription of both states and controls by Finite Elements in Time (DFET)[7]. A set of additional parameters, not included among states and controls, are allowed and can be used for a combined optimization of both the trajectory and other quantities peculiar to the original optimal control problem (parametric optimization). In particular, in this paper, the orbital elements of each hyperbola are treated as additional parameters and, opposite to the work of Betts, swing-by trajectories are not transcribed with collocation but using multiple shooting.

The method is successfully applied to the design of a mission to planet Mercury exploiting one or more swing-bys of Venus and Mercury itself. Different options to reach planet Mercury are analyzed, each characterized either by a different sequence of swing-bys or by a different launch opportunity. In particular, two strategies involving two swing-bys of Venus and two or three swing-bys of Mercury have been studied. In order to validate the proposed optimization approach, a comparison is made between the solution obtained with DFET and the result coming from optimal control theory for a reduced problem.

**PROBLEM FORMULATION**

The problem is modeled in two different ways of increasing complexity. First as a reduced two-body problem, with the Sun as primary and the swing-bys treated as singular events, instantaneous and with no variation in position. Then as a full three dimensional problem with swing-bys treated as actual three-dimensional trajectories in space and time including perturbations. The former solution is used to provide a first guess to the latter.

The date and the position of the encounter, are completely free, as well as the departure date from the Earth and the arrival at Mercury. The only piece of information that must be provided is the number and name of celestial bodies used for the gravity maneuvers. The sequence and type of celestial bodies employed distinguishes each different strategy to reach Mercury. Although guessing the swing-bys bodies could be regarded as a limitation, from a mission design point of view, it allows the analyst to design each swing-by in the most appropriate way, inserting even special conditions (e.g. coast arcs, before each encounter, required for navigation), since the early design stages.

In order to take into account swing-bys, the trajectory has been split into several phases, each phase corresponding to a trajectory arc connecting two planets. On each phase a particular collocation technique based on Finite Elements in Time has been used to transcribe differential equations, governing the dynamics of the spacecraft, into a set of algebraic nonlinear equations and to parameterize controls. When treating swing-bys as full three-dimensional trajectories, a local reference frame is taken to describe the gravity assist maneuvers. Incoming conditions, at the sphere of influence, represent final conditions for the phase preceding the swing-by and outgoing conditions, at the sphere of influence, represent initial conditions for the subsequent phase. Within the sphere of influence, hyperbolas are propagated backward and forward in time from the pericenter in a local reference frame taking into account perturbations from the Sun. In this way collocation and multiple shooting are combined in a unique approach reducing the number of collocation points required but retaining robustness.

All the phases are then assembled together, forming a single NLP problem. Each phase is linked to the preceding one and to the following one by the appropriate set of boundary conditions computed by the relative swing-by trajectory. The resulting nonlinear programming problem (NLP) is highly sparse and can be solved efficiently by any sparse sequential programming algorithm. In the present case, the SNOPT[8] library has been used.

In the following paragraphs the dynamic model used to describe the trajectory between two encounters and the two different swing-by models employed are presented.

**Dynamics**

A spacecraft is modeled as a point mass subject to the gravity attraction of the Sun and to the thrust provided by one or more low-thrust engines. Introducing the position vector $\mathbf{r}=\{r_x, r_y, r_z\}^T$, the velocity vector $\mathbf{v}=\{v_x, v_y, v_z\}^T$ and the propellant mass $m_p$, the state and the control vectors are then defined as follows:

$$\begin{aligned} \mathbf{x} &= \{r_x, r_y, r_z, v_x, v_y, v_z, m_p\}^T; \\ \mathbf{u} &= \{u_x, u_y, u_z\}^T \end{aligned} \quad (1)$$

The motion of the spacecraft is described in the J2000 reference frame centered in the Sun (Figure1). The three components of the thrust vector $\mathbf{u}$ represent the control:

$$\begin{aligned} \dot{\mathbf{r}} &= \mathbf{v} \\ \dot{\mathbf{v}} &= \nabla U(\mathbf{r}) + \frac{\mathbf{u}}{m_d + m_p} \end{aligned} \quad (2)$$

$$\dot{m}_p = -\frac{u}{I_{sp}g_0}$$

where $I_{sp}$ is the specific impulse of the engine and $g_0$ the gravity constant on Earth surface, and the gravity potential of the Sun is a function of the position vector $\mathbf{r}$ and the gravity constant $\mu$:

$$U(\mathbf{r}) = \frac{\mu}{|\mathbf{r}|} \quad (3)$$

The mass of the spacecraft is divided into propellant mass $m_p$ and dry mass $m_d$. An upper bound $T_{max}$ and a lower bound $T_{min}$ was put on the thrust magnitude:

$$T_{min} \leq u = \sqrt{u_x^2 + u_y^2 + u_z^2} \leq T_{max} \quad (4)$$

The upper bound is the maximum level of thrust provided by the selected low-thrust engine, the lower was taken $10^{-4}$ times $T_{max}$ to avoid singularities in the Hessian matrix when minimum mass problems are solved. The control vector $\mathbf{u}$ can be better represented in a local reference frame centered in the spacecraft by decomposing it into a tangential component $u_v$ aligned with the velocity vector, a normal component u_n, normal to the trajectory and a bi-normal component u_h, normal to the orbital plane. In this reference frame the elevation angle ϕ is defined as the angle between the control vector **u** and the plane tangential to the trajectory containing u_v and u_h. In the same reference, the azimuth angle α is defined as the angle between the projection of the control vector in the tangent plane and the velocity vector **v** (see Figure 1).

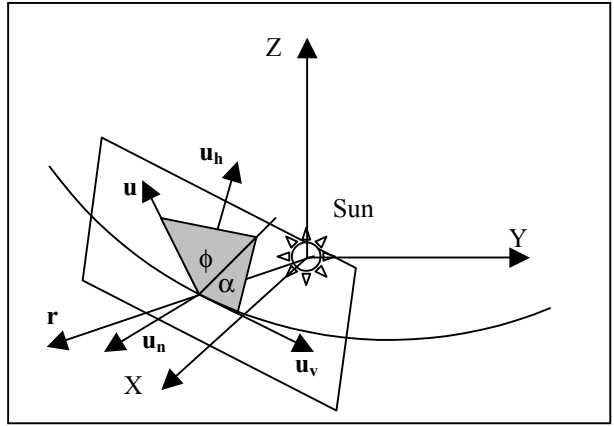

Figure 1. Inertial reference frame centered in the Sun: the xy plane is the ecliptic plane and x-axis points toward the 2000 mean vernal equinox.

**Swing-by**

The simplest way to model a gravity assist maneuver is to resort to link-conic approximation: the sphere of influence of a planet is assumed to have zero radius and the gravity maneuver is considered instantaneous. Therefore the instantaneous position vector is not affected by the swing-by:

$$\mathbf{r}_i = \mathbf{r}_o = \mathbf{r}_p \qquad (5)$$

where $\mathbf{r}_i$ is the incoming heliocentric position, $\mathbf{r}_o$ is the outgoing heliocentric position vector and $\mathbf{r}_p$ is the planet position vector, all taken at the epoch of the encounter. For an ideal hyperbolic orbit, not subject to perturbations or Δv maneuvers, the modulus of the incoming relative velocity must be equal to the modulus of the outgoing relative velocity:

$$\tilde{v}_i = \tilde{v}_o \qquad (6)$$

Furthermore the outgoing relative velocity vector is rotated, due to gravity, of an angle 180-2β with respect to the incoming velocity vector and therefore the following relation must hold:

$$\tilde{\mathbf{v}}_o^T \tilde{\mathbf{v}}_i = -\cos(2\beta)\tilde{v}_i^2 \qquad (7)$$

where the complementary angle of rotation of the velocity is defined as:

$$\beta = acos\left(\frac{\mu}{\tilde{v}_i^2 \tilde{r}_p + \mu}\right) \quad (8)$$

All quantities with a tilde are relative to the swing-by planet and $\tilde{r}_p$ is the periapsis radius of the swing-by hyperbola (see Figure 2).

**Numerical Propagation**

After a solution has been obtained with the link-conic model, a second solution is computed substituting the simple link-conic approximation with a fully 3d numerical propagation of the swing-by hyperbolas. Each swing-by is treated as a new phase that has to be linked to the incoming part of the trajectory and to the outgoing part of the trajectory at the sphere of influence. Swing-bys are not propelled and therefore there is no need to introduce a control on the thrust vector along the swing-by hyperbola. Thus two reference frames are used and two dynamical models. The first one is a heliocentric reference frame and the spacecraft is subject to the gravity attraction of the Sun and to the thrust of the SEP engine. The second is centered into the swing-by planet and the spacecraft is subject to the gravity attraction of the swing-by planet and to third body perturbations coming from the Sun. Thus the dynamics of the spacecraft within the sphere of influence is governed by the following differential equation:

$$\frac{d\tilde{\mathbf{x}}}{dt} = \tilde{\mathbf{F}}(\tilde{\mathbf{x}},t) = \begin{cases} \tilde{\mathbf{v}} \\ -\frac{\mu_P}{\tilde{r}^3}\tilde{\mathbf{r}} - \mu_S\left(\frac{\mathbf{d}}{d^3} + \frac{\mathbf{r}_S}{r_S^3}\right) \end{cases} \quad (9)$$

where $\mathbf{d}$ is the spacecraft-Sun vector and $\mathbf{r}_S$ is the position vector of the Sun in the planetocentric reference frame. In order to increase robustness, orbital parameters for each hyperbola are not derived from incoming conditions but are included into the set of NLP parameters and then optimized. Hyperbolas are propagated backward in time from the pericenter up to the sphere of influence, where they are linked to the incoming trajectory, and forward in time up to the sphere of influence, linked to the outgoing trajectory. The values of the orbital parameters are then optimized in order to satisfy matching conditions on the sphere of influence. A first guess value for the parameter is obtained from the previous solution, the semimajor axis and the eccentricity can be easily derived from the incoming velocity modulus and from the pericenter radius:

$$a = -\frac{\mu}{\tilde{v}_i^2}; \quad e = 1 - \frac{\tilde{r}_p}{a} \quad (10)$$

The incoming and the outgoing velocity vectors must both lie in the orbital plane and therefore:

$$\mathbf{h} = \tilde{\mathbf{v}}_i \wedge \tilde{\mathbf{v}}_o \ ; \ \mathbf{N} = \left[-h_y/h, h_x/h\right]^T$$
$$i = \mathrm{acos}\frac{h_z}{h}; \quad \Omega = \mathrm{atan}\frac{h_x}{-h_y} \quad (11)$$

the apsidal line $[b_x, b_y, b_z]$ must bisect the angle between the incoming and the opposite of the outgoing vector and must lie in the orbital plane, therefore the following linear system must hold:

$$\begin{bmatrix} h_x & h_y & h_z \\ \tilde{v}_x^i & \tilde{v}_y^i & \tilde{v}_z^i \\ \tilde{v}_x^o & \tilde{v}_x^o & \tilde{v}_x^o \end{bmatrix} \begin{Bmatrix} b_x \\ b_y \\ b_z \end{Bmatrix} = \begin{bmatrix} 0 \\ \cos\beta \\ -\cos\beta \end{bmatrix} \quad (12)$$

The three components of the apsidal axis are obtained solving the previous linear problem while the anomaly of the pericenter can be computed as the angular distance between the apsidal line and the line of the nodes:

$$\omega = \mathrm{acos}\frac{\mathbf{b} \cdot \mathbf{N}}{\|b\|\|N\|} \quad (13)$$

In addition to the five orbital parameters, for each hyperbola the time spent within the sphere of influence is derived from the semimajor axis and the eccentricity:

$$\cosh H = \frac{1}{e}\left(1 - \frac{\tilde{r}_i}{a}\right) \quad (14)$$

$$\Delta t = (e\sinh H - H)\sqrt{\frac{a^3}{\mu}} \quad (15)$$

This value is used to integrate backward in time the state vector computed at the pericenter of the hyperbola up to the sphere of influence and forward in time the same state vector up to the sphere of influence. The state vector at the pericenter of the hyperbola is computed from the orbital parameters:

$$[\tilde{\mathbf{r}}_p, \tilde{\mathbf{v}}_p]^T = f(a, e, i, \omega, \Omega, 0) \quad (16)$$

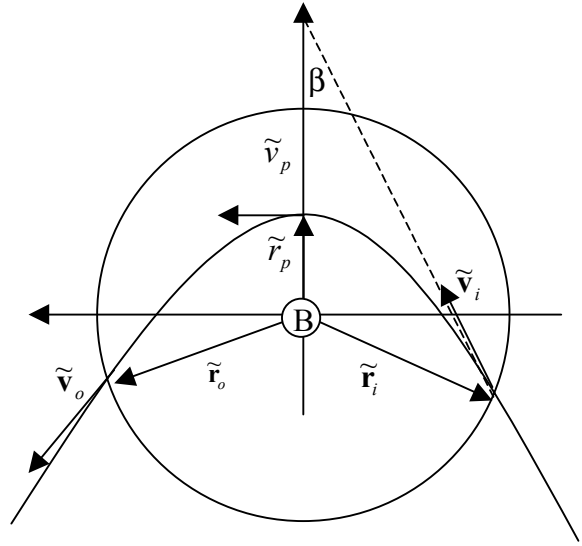

Figure 2. Swing-by model and reference frame

Therefore at the sphere of influence of a body B, with state vector $[\mathbf{r}_B, \mathbf{v}_B]^T$, the following set of matching constraints must be satisfied:

$$\begin{aligned}\tilde{\mathbf{v}}_i &= \mathbf{v}_i - \mathbf{v}_B(t-\Delta t) \\ \tilde{\mathbf{r}}_i &= \mathbf{r}_i - \mathbf{r}_B(t-\Delta t)\end{aligned} \quad (17)$$

$$\begin{aligned}\tilde{\mathbf{v}}_o &= \mathbf{v}_o - \mathbf{v}_B(t+\Delta t) \\ \tilde{\mathbf{r}}_o &= \mathbf{r}_o - \mathbf{r}_B(t+\Delta t)\end{aligned} \quad (18)$$

where incoming and outgoing relative position and velocity vectors are obtained integrating respectively from $t$ to $t-\Delta t$ and from $t$ to $t+\Delta t$ the differential equations:

$$\begin{aligned}\tilde{\mathbf{x}}_i &= \tilde{\mathbf{x}}_p + \int_t^{t-\Delta t} \tilde{\mathbf{F}}(\tilde{\mathbf{x}},t)\,dt \\ \tilde{\mathbf{x}}_o &= \tilde{\mathbf{x}}_p + \int_t^{t+\Delta t} \tilde{\mathbf{F}}(\tilde{\mathbf{x}},t)\,dt\end{aligned} \quad (19)$$

Figure 2 reports a sketch of the model adopted for swing-bys.

**Thrust Model**

The thrust provided by the engine is determined taking into account the specific thrust $F_{sp}$, the effective input power $P_{in}$ provided by the power system and an efficiency coefficient $\eta_e$:

$$F_{max} = \eta_e P_{in} F_{sp} \quad (20)$$

The effective input power is given by the effective power produced by the solar arrays minus the power required by the spacecraft $P_{ss}$:

$$P^*_{in} = P_{eff} - P_{SS} \tag{21}$$

In order to take into account the degradation of the solar arrays due to temperature and the reduced power due to the increasing distance from the Sun, the power provided by the solar arrays during the transfer trajectory is here expressed in the approximate form:

$$P_{eff} = \eta_S \frac{P_{1AU}}{R_S^2}[1 - C_T(T_S - T_0)]\cos\alpha \tag{22}$$

where $P_{IAU}$ is the power at one Astronomical Unit, $T_s$ is the temperature of solar arrays, $R_S$ is the distance from the Sun, $T_0$ the reference temperature, $C_T$ is the temperature coefficient which expresses the reduced performance of the panel with temperature increase, $\eta_S$ is a coefficient to account for all other degradations sources and $\alpha$ is the solar array Sun aspect angle, i.e. the angle between the normal to the cell surface and the Sun direction. The steady state surface temperature of the solar panels is here taken as function of the distance from the Sun assuming purely radiative heat transfer for the solar array:

$$T_S = \left[\frac{S_0 \alpha_s \cos\alpha}{R_S^2 \sigma \kappa \varepsilon}\right]^{0.25} \tag{23}$$

where $S_0$ is the solar constant at 1 AU, $\sigma$ is the Stefan-Boltzmann constant, $\alpha_s$ is the surface absorbivity is the solar spectrum ad $\varepsilon$ is the surface emissivity is the infrared spectrum, $\kappa$ is a coefficient which takes into account the surface area radiating in the infrared spectrum, with respect to the one that receives the solar input. A maximum power that can be handled by the PPU is assumed to represent the upper limit for the engine thrust.

$$P_{in} = min(P^*_{in}, P_{max}) \tag{24}$$

The required power is dimensioning for the design of the solar arrays and power system and therefore it provides estimation for the overall dry mass of the spacecraft.

**OPTIMISATION APPROACH**

A general trajectory design problem can be decomposed in $M$ phases, each one characterized by a time domain $D^j$, with $j=1,...,M$, a set of $m$ dynamic variables **x**, a set of $n$ control variables **u** and a set of $l$ parameters **p**. Furthermore, each phase $j$ may have an objective function

$$J^j = \phi^j(\mathbf{x}_0^b, \mathbf{x}_f^b, t_f, \mathbf{p}) + \int_{t_i}^{t_f} L^j(\mathbf{x}, \mathbf{u}, \mathbf{p}) dt \tag{25}$$

a set of dynamic equations

$$\dot{\mathbf{x}} - \mathbf{F}^j(\mathbf{x}, \mathbf{u}, \mathbf{p}, t) = 0 \tag{26}$$

a set of algebraic constraints on states and controls

$$\mathbf{G}^j(\mathbf{x}, \mathbf{u}, \mathbf{p}, t) \geq 0 \tag{27}$$

and a set of boundary constraints

$$\psi^j(\mathbf{x}_0^b, \mathbf{x}_f^b, \mathbf{p}, t)\Big|_{t_0}^{t_f} \geq 0 \tag{28}$$

Among boundary constraints a set of inter-phase link constraints exist that are used to assemble all phases together

$$\psi^j(\mathbf{x}_j^b, \mathbf{x}_{j-1}^b, \mathbf{p}, t) \geq 0 \tag{29}$$

The time domain $D(t_0, t_f) \subset \Re$ relative to each phase $j$ can be further decomposed into $N$ finite time elements $D^j = \bigcup_{i=1}^{N} D_i^j(t_{i-1}, t_i)$ and, on each time element $D_i^j$, states and controls $[\mathbf{x}, \mathbf{u}]$ can be parameterized as follows:

$$\begin{Bmatrix} \mathbf{x} \\ \mathbf{u} \end{Bmatrix} = \sum_{s=1}^{p} f_s(t) \begin{Bmatrix} \mathbf{x}_s \\ \mathbf{u}_s \end{Bmatrix} \tag{30}$$

where the basis functions $f_s$ are chosen within the space of polynomials of order $p-1$:

$$f_s \in P^{p-1}(D_i^j) \tag{31}$$

Therefore in general a finite element is defined by a sub-domain $D_i^j$, and by a sub-set of parameters $[\mathbf{x}_s, \mathbf{u}_s, \mathbf{p}]$. A group of finite elements forms a phase and a group of phases forms the original optimization problem. Notice that additional parameters $\mathbf{p}$ may appear in all constraint equations depending on their meaning in the optimization problem. Furthermore, it should be noticed that each phase could be grouped in sequence or in parallel with the other phases depending on its time domain and on the inter-phase link constraints that pass information among phases. Thus two phases can share the same time domain but have different parameterizations.

Now taking a general phase, in order to integrate differential constraints (26), on each finite element $i$, differential equations are transcribed into a weighted residual form considering boundary conditions of the weak type:

$$\int_{t_i}^{t_{i+1}} \left\{ \dot{\mathbf{w}}^T \mathbf{x} + \mathbf{w}^T \mathbf{F}^j \right\} dt - \mathbf{w}_{i+1}^T \mathbf{x}_{i+1}^b + \mathbf{w}_i^T \mathbf{x}_i^b = 0 \qquad i = 1, \ldots, N-1 \qquad (32)$$

where **w**(t) are generalized weight (or test) functions defined as:

$$\mathbf{w} = \sum_{s=1}^{p+1} g_s(t) \mathbf{w}_s \qquad (33)$$

where $g_s$ are taken within the space of polynomials of order $p$:

$$g_s \in P^p(D_i^j) \qquad (34)$$

Now the problem is to find the vector $\mathbf{x}_s \in \Re^{p*m}$, the vector $\mathbf{u}_s \in \Re^{p*n}$, the vector $\mathbf{p} \in \Re^l$ and $\mathbf{x}_f^b$ and $\mathbf{x}_0^b$ $\in \Re^m$ that satisfy variational equation (32) along with algebraic and boundary constraints:

$$\mathbf{G}^j(\mathbf{x}, \mathbf{u}, \mathbf{p}, t) \geq 0 \qquad (35)$$

$$\psi^j(\mathbf{x}_0^b, \mathbf{x}_f^b, \mathbf{p}, t)\Big|_{t_0}^{t_f} = 0 \qquad (36)$$

where quantities $\mathbf{x}_s$, and $\mathbf{u}_s$ are called internal node values, while $\mathbf{x}_f^b$, $\mathbf{x}_0^b$ are called boundary values. Notice that generally the order $p$ of the polynomials can be different for states and controls. In a more general way the domain $D^j$ could be decomposed as a union of smooth images of the reference time interval [-1,1] where a reference parameter $\tau$ is defined as:

$$\tau = 2 \frac{t - t_{i-1/2}}{t_i - t_{i-1}} = 2 \frac{t - t_{i-1/2}}{\Delta t_i} \qquad (37)$$

Polynomials $f_s$ and $g_s$ are constructed using Lagrangian interpolants associated with internal Gauss-type nodes. Generally speaking if $\{\xi_s\}^p_{s=1}$ are the set of Gauss points on the reference interval [-1,1], $f_s(\tau)$ will be the Lagrangian interpolating polynomial vanishing at all Gauss points except at $\xi_s$ where it equals one.

Each integral of the continuous forms (25) and (32) is then replaced by a q-points Gauss quadrature sum, where q is taken equal to p. Therefore the objective function (25) becomes a sum of $N$ Gauss quadrature formulas:

$$J^j = \phi^j(\mathbf{x}_0^b, \mathbf{x}_f^b, t_f) + \sum_{i=1}^{N} \sum_{k=1}^{q} \sigma_k L_k^j \frac{\Delta t_i}{2} \qquad (38)$$

while integral (32) is split into $N$ integrals of the form:

$$\sum_{k=1}^{q} \sigma_k \left[ \dot{\mathbf{w}}_k(\tau_k)^T \mathbf{x}(\tau_k) + \mathbf{w}_k(\tau_k)^T \mathbf{F}_k^j \frac{\Delta t_i}{2} \right] - \mathbf{w}_{p+1}^T \mathbf{x}_{i+1}^b + \mathbf{w}_1^T \mathbf{x}_i^b = 0 \qquad i = 1, \ldots, N-1 \qquad (39)$$

where $\sigma_k$ are Gauss weights and parameters $\mathbf{x}^b_{i-1}$ and $\mathbf{x}^b_i$ are boundary values at the beginning and at the end of each element. For sake of simplicity, the following notation has been introduced:

$$L^j_k = L^j(\mathbf{x}_s f_s(\tau_k), \mathbf{u}_s f_s(\tau_k), \mathbf{p}, \tau_k); \quad \mathbf{F}^j_k = \mathbf{F}^j(\mathbf{x}_s f_s(\tau_k), \mathbf{u}_s f_s(\tau_k), \mathbf{p}, \tau_k) \quad (40)$$

Here controls are parameterized using the same set of points used for integration while states are always collocated on Gauss-Lobatto nodes. Numerical quadrature of the integral Eq. (32) and integral (25) can be then performed either by Gauss Lobatto rule or by Gauss-Legendre rule. The former choice of quadrature formulas collocates controls on the same set of nodes as states while the latter collocates controls on a different set. The advantage of the latter is the higher integration order, which allows a lower number of collocation nodes. Whatever $f_s$ and $g_s$ are, the linear part of Eq. (39) can be always integrated only once before the optimization process begins. Now Eq. (39) must be satisfied for every arbitrary value of virtual quantity $\mathbf{w}_k$, as a consequence each element equation is developed into $p+1$ equations:

$$\sum_{k=1}^{q} \sigma_k \mathbf{F}^j_k \frac{\Delta t_i}{2} \begin{Bmatrix} g_1(\tau_k) \\ \vdots \\ g_{p+1}(\tau_k) \end{Bmatrix} + \begin{bmatrix} \sum_{k=1}^{q} \sigma_k \dot{g}_1(\tau_k) f_1(\tau_k) & \cdots & \sum_{k=1}^{q} \sigma_k \dot{g}_1(\tau_k) f_p(\tau_k) \\ \vdots & \ddots & \vdots \\ \sum_{k=1}^{q} \sigma_k \dot{g}_{p+1}(\tau_k) f(\tau_k)_1 & \cdots & \sum_{k=1}^{q} \sigma_k \dot{g}_{p+1}(\tau_k) f_p(\tau_k) \end{bmatrix} \begin{Bmatrix} \mathbf{x}_1 \\ \vdots \\ \mathbf{x}_p \end{Bmatrix} = \begin{bmatrix} -\mathbf{x}^b_i \\ 0 \\ \mathbf{x}^b_{i+1} \end{bmatrix} \quad (41)$$

System of Eqs. (41) is written for each element, all the elements are then assembled matching the final boundary node of one element to the initial one of the next element. For continuous solution, in order to preserve the continuity of the states, at matching points, the following condition must hold:

$$\mathbf{x}^b_i = \mathbf{x}^b_{i+1} \qquad i=1,\ldots,N\text{-}2 \quad (42)$$

Thus all the boundary quantities (42) cancel one another except for those at the initial and final times. Algebraic constraint equation (35) can be collocated directly at Gauss nodal points:

$$\mathbf{G}^j_s(\mathbf{x}_s(\xi_s), \mathbf{u}_s(\xi_s), \mathbf{p}, \xi_s) \geq 0 \quad (43)$$

The resulting set of non-linear algebraic equations, assembling all the phases, along with discretised objective function (38) can be seen as a general non-linear programming problem (NLP) of the form:

$$\min \ J(\mathbf{y}) \quad (44)$$

subject to

$$\begin{aligned} \mathbf{c}(\mathbf{y}) &\geq 0 \\ \mathbf{b}_l &\leq \mathbf{y} \leq \mathbf{b}_u \end{aligned} \quad (45)$$

where, $\mathbf{y}$ is the vector of NLP variables, $J(\mathbf{y})$ is the objective function to be minimized, $\mathbf{c}(\mathbf{y})$ is a vector of non-linear constraints and $\mathbf{b}_l$ and $\mathbf{b}_u$ are respectively lower and upper bounds on NLP variables. The $N*(p+1)*n$ algebraic Eqs. (41) taken for each phase, along with system (43), represent the $\mathbf{c}(\mathbf{y})$ constraint of the nonlinear problem while $\mathbf{y}=[\mathbf{x}_s,\mathbf{u}_s,\mathbf{x}^b_0,\mathbf{x}^b_f,t_0,t_f,\mathbf{p}]$ the NLP variables. Notice that the present formulation is discontinuous because continuity at boundaries of each element is only weakly enforced. This means that, generally, there is a jump between the internal nodes and the boundary nodes. This allows the control, for which no continuity requirement is imposed, to be discontinuous at boundaries.

## RESULTS

In the following DFET, with multiphase and parametric optimization approach, has been used to design different trajectories toward Mercury. Two swing-bys of Venus and two swing-bys of Mercury (EVVMM strategy) characterize the first trajectory. The second is characterized by two swing-bys of Venus and a sequence of three swing-bys of Mercury, resonant with the motion of the planet (EVVMMM strategy).

In the former case two options for the launch have been taken into consideration: a Soyuz launch with Fregat upper stage and an Ariane 5 launch with Vinci upper stage. For the latter case only Ariane 5 has been considered as possible launcher. For all cases, power supply characteristics have been kept as summarized in table 1.

Table 1. Power system characteristics

| Parameter | Value |
|---|---|
| $\eta_e$ | 0.9 |
| $\eta_S$ | 1 |
| $P_{IAU}$ | 11.2 kW |
| $C_T$ | 3*10-4 K-1 |
| $T_0$ | 290 K |
| $\kappa$ | 1.3 |
| $\varepsilon$ | 0.86 |
| $\alpha$ | 0.86 |
| $T_{max}$ | 423 K |
| $P_{max}$ | 15 kW |
| $P_{SS}$ | 300 W |

**Test Case: Comparison Between DFET and Indirect Solution**

At first a comparison has been made with a solution computed using Pontryagin maximum principle[5] (indirect approach). The case studied aims to inject a spacecraft, with an initial wet mass of 1280 kg, into a high elliptical orbit around Mercury starting from a fixed launch date and maximizing the arrival mass. A Soyuz launcher with Fregat upper stage provides 2.3 km/s of infinite velocity at departure from the Earth. The spacecraft is equipped with a SEP module which provides a maximum nominal thrust of 0.34N with an Isp of 3200 s. The arrival asymptotic velocity is required to be 360 m/s. The strategy employed consists of two swing-bys of Venus in order to reduce the perihelion, then after four complete revolutions around the Sun a first encounter with Mercury reduces the aphelion, while a second encounter with Mercury injects the spacecraft in the right orbit leading to the desired final encounter with Mercury. Therefore the problem is split into five phases with four swing-bys. A first guess for the first phase, from the Earth to Venus, was computed solving the corresponding Lambert problem. The second phase, from Venus to Venus, was guessed simply propagating forward in time the orbit of Venus, while the third phase, from Venus to Mercury, was generated at first integrating forward in time with a constant thrust opposite to the velocity vector. After that a first solution was computed without any upper bound on the thrust and using as objective function the square of the modulus of the control vector.

This intermediate solution satisfies all the boundary conditions keeping the average value of the thrust quite below the upper limit. A second solution is then computed minimizing the final mass imposing the required upper limit on the thrust modulus. This results into a bang-bang switching structure with 12 thrust arcs and 22 switching points. Transfer trajectory obtained with DFET is represented in Figures 3 and 4 along with the indirect solution.

The time history of the velocity component perpendicular to the ecliptic plane is represented in Figure 5, while in Figures 6 and 7 the thrust elevation and magnitude time histories are reported respectively. In these two last figures the complexity of the switching structure can be clearly seen. It should be noted that no a priori information about the switching structure were provided to DFET that was able to reconstruct correctly the sequence of thrust and coast arcs. Compared to the indirect method, DFET solution presents an additional thrust arc, but with a reduced thrust of 0.01N. This additional arc is probably due to a lack of accuracy in the exact reconstruction of the other thrust arcs. A better fitting of the switching structure may solve the problem. A further difference between the two solutions is in the arrival conditions. As can be noticed in the velocity plot, there is a small difference in final values. However the arrival date is the same for both solutions, as can be read in Table 2, where a comparison between the two solutions is reported. This difference may be due to a small difference in the ephemeris data for Mercury. In the same table, the encounter date for each swing-by is reported showing a good matching of the two solutions, with a difference of just few ours, less than 0.3 days.

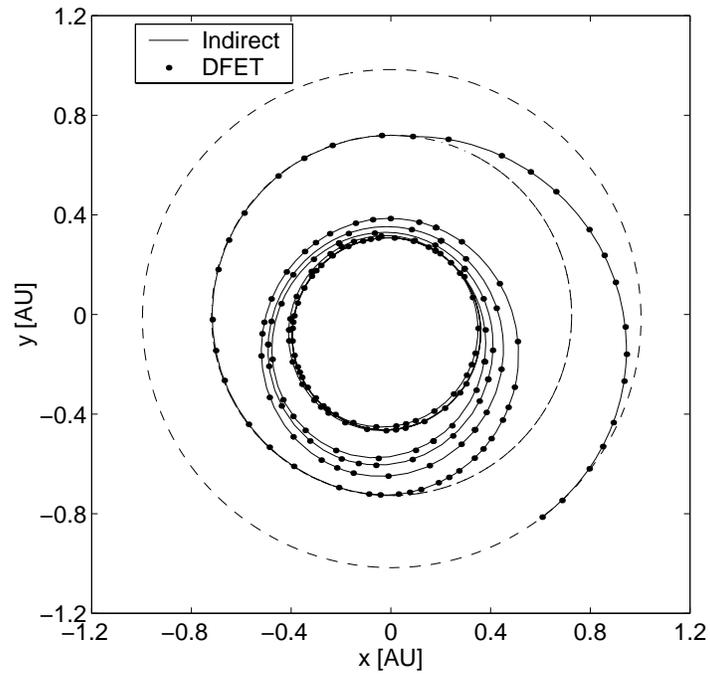

Figure 3. Comparison between DFET (dotted line) and indirect method (solid line): ecliptic plane.

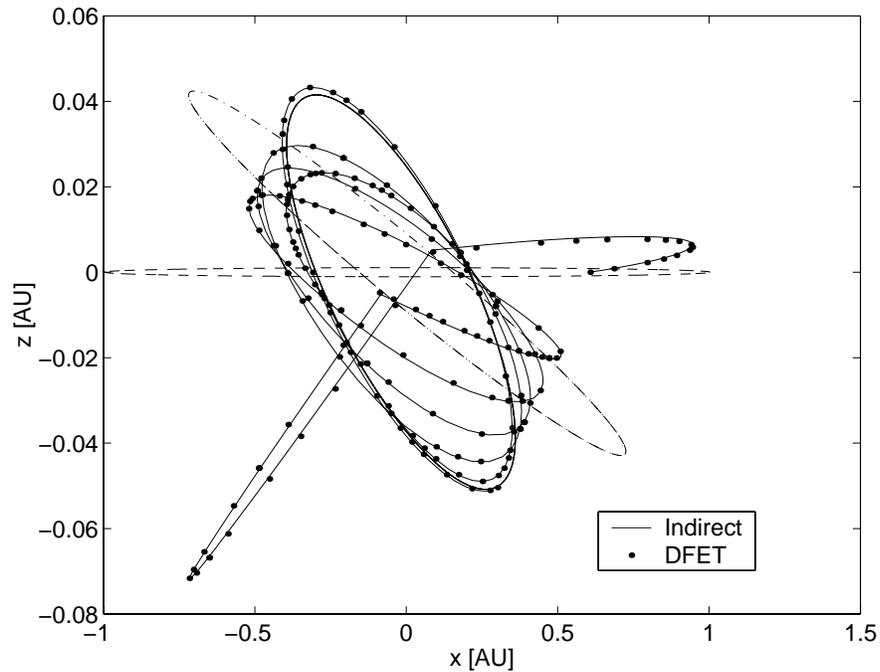

Figure 4. Comparison between DFET (dotted line) and indirect method (solid line): plane orthogonal to ecliptic.

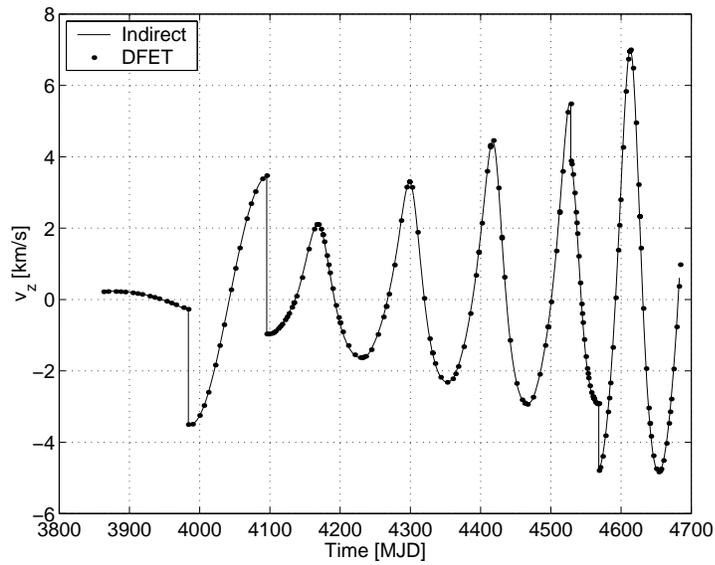

Figure 5. Comparison between DFET (dotted line) and indirect method (solid line): $v_z$ axis.

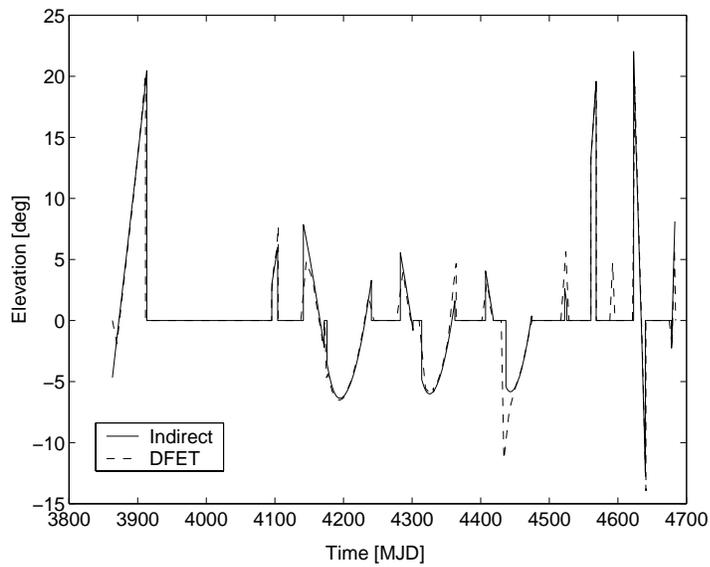

Figure 6. Comparison between DFET (dotted line) and indirect method (solid line): thrust elevation.

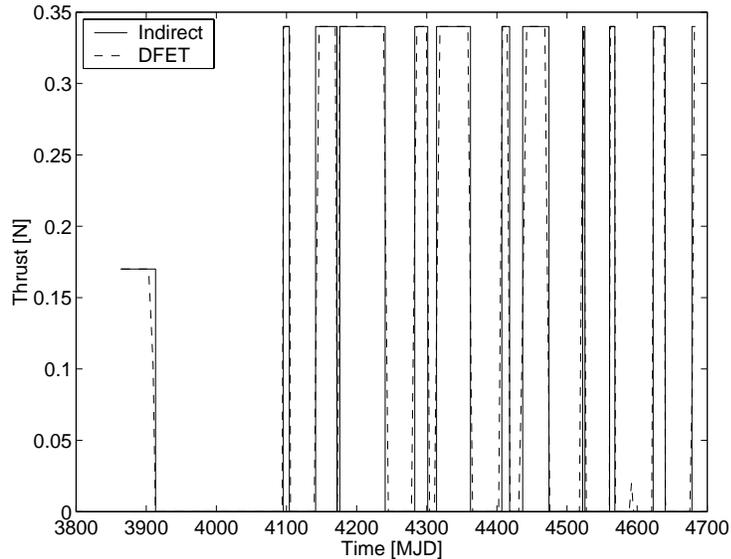

Figure 7. Comparison between DFET (dotted line) and indirect method (solid line): thrust magnitude.
16

It should be noted that the need to provide a first guess solution is due to the expected relevant number of minima and not to convergence difficulties. In fact it has been verified that even though DFET approach converges to a solution even without a particularly good initial guess, demonstrating to be quite robust, never the less it might converge to different solutions depending on the initial guess provided.

Table 2. Comparison between DFET and indirect solution

| DATA | DFET | | INDIRECT | |
|---|---|---|---|---|
| $v_\infty$ | 2.3 km/s | | 2.3 km/s | |
| Max Thrust | 0.17-0.34 | | 0.17-0.34 | |
| Isp | 3200 s | | 3200 s | |
| Launch Date | 30 July 2010 | | 30 July 2010 | |
| Initial Mass | 1280 kg | | 1280 kg | |
| Final Mass | 1021.1 kg | | 1020.9 kg | |
| Payload Mass | 584.52 kg | | 584.34 kg | |
| Arrival Velocity | 0.360 km/s | | 0.360 km/s | |
| Arrival Date | 27 Oct 2012 | | 27 Oct 2012 | |
| | Altitude | Date (MJD) | Altitude | Date (MJD) |
| Venus Swing-by | 2674.9 km | 3983.7 | 2761.3 km | 3983.9 |
| Venus Swing-by | 300 km | 4095.3 | 300 km | 4095.5 |
| Mercury Swing-by | 200 km | 4528.6 | 200 km | 4528.5 |
| Mercury Swing-by | 200 km | 4568.5 | 200 km | 4568.2 |

**Application Example 1: Soyuz case EVVMM strategy**

Once the DFET approach has been tested reproducing correctly the solution obtained from control theory, the actual problem of a transfer from the Earth to Mercury has been solved introducing a full model for the swing-bys and taking into account even performances of the launcher at departure. In addition it has been considered that before each encounter a coast arc of about 30 days should be inserted to allow a good orbit determination especially before each swing-by. Once the probe has reached Mercury, the SEP module is jettisoned and the final break is performed with a chemical engine characterized by an Isp of 317 s and a maximum thrust of 400 N. Therefore the objective function is the payload mass that should be delivered into a 12000x400 km orbit around Mercury after the maneuver with chemical propulsion. A model for Soyuz performances has been introduced as boundary constraint on the departure asymptotic velocity and mass while a model for the final break has been introduced to compute the objective function as a function of the arrival velocity. The resulting trajectory is represented in Figure 8, the solution obtained with DFET has been propagated forward in time using a Runge-Kutta 5th order variable step-size integrator to verify the quality of the solution. The DFET solution is propagated in a heliocentric reference frame with an n-body gravity model (i.e. including actual gravity of each planet). The solid line represents thrust arcs while the dashed line represents coast arcs. As can be clearly seen the imposed period of 30 days before each encounter is satisfied. A circle marker represent the departure from the Earth while stars represent swing-bys and a diamond marker represents the arrival at Mercury. Two star markers represent respectively the entry into the sphere of influence end the exit point from the sphere of influence for each swing-by. The time elapsed within the sphere of influence varies from about 2.5 days for Venus to 1.2 days for Mercury. Each swing by is fully numerical, the result obtained after propagation is represented in Figure 8 for the entire trajectory in the ecliptic plane and in Figure 9 where the time history of velocity component perpendicular to the ecliptic plane is represented. The effects of the four swing-bys is quite evident and, as an example, a close up of the second encounter with Venus is represented in Figure 10, showing the accuracy of the DFET solution. In fact the error at the sphere of influence where the propagated hyperbola are linked to the transfer arcs is less than 1e-3.

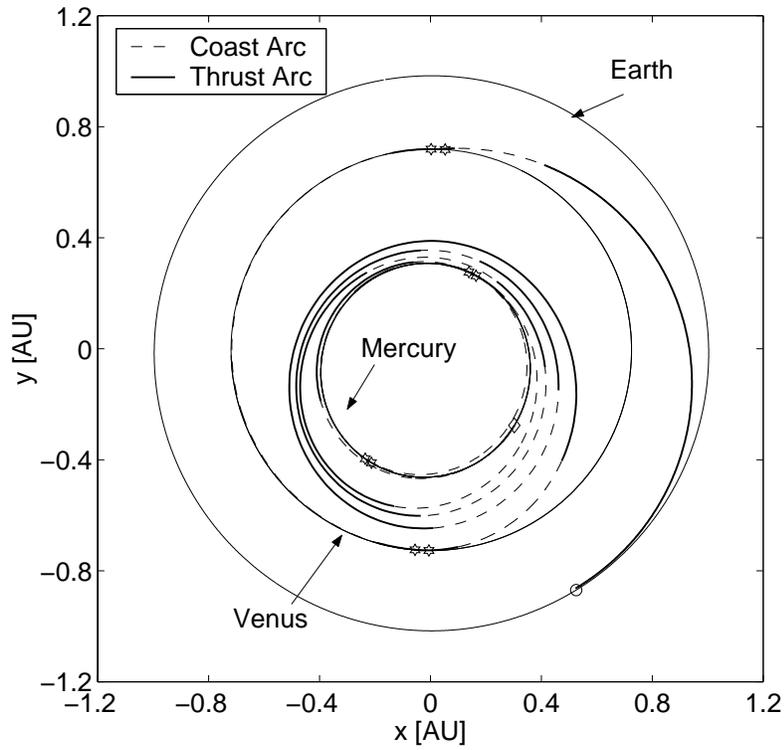

Figure 8. Trajectory in the ecliptic plane for the Soyuz case with EVVMM strategy

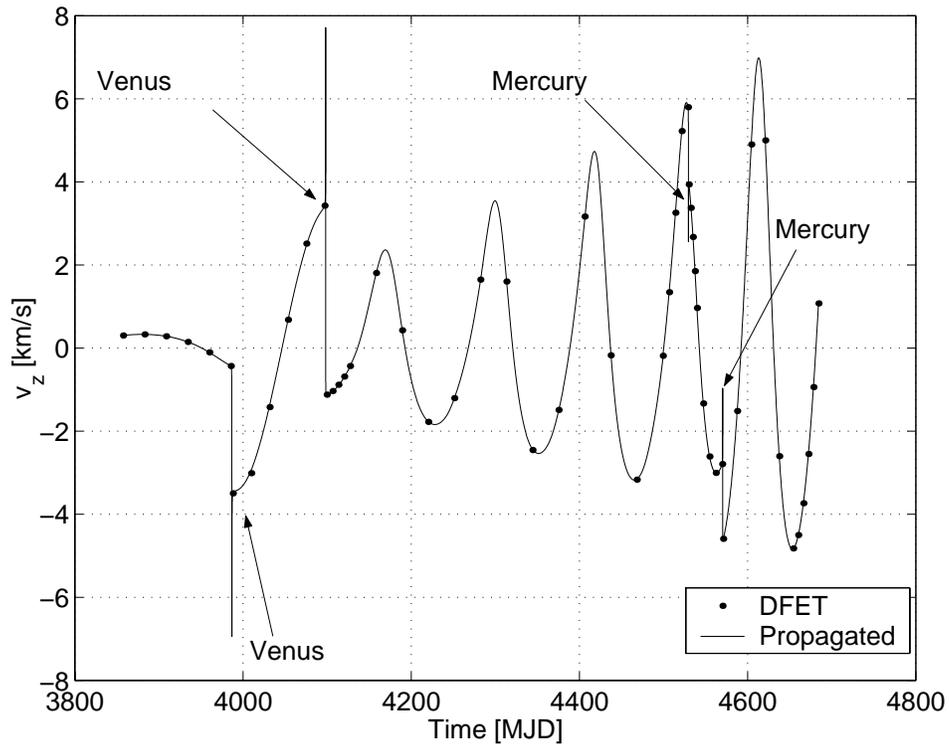

Figure 9. Time history of the velocity component along z-axis

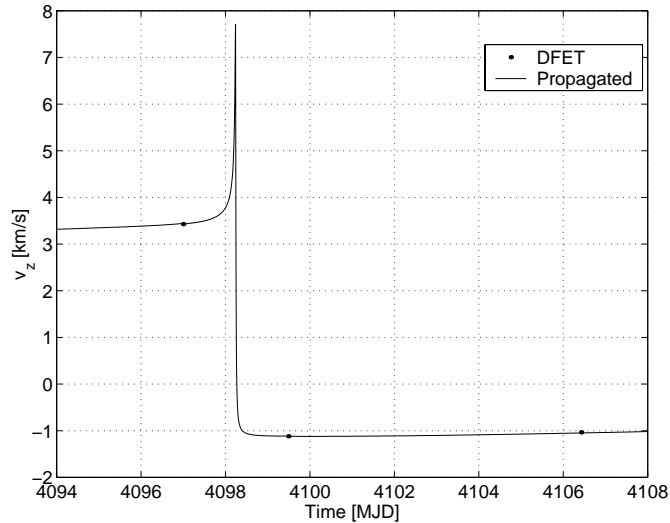
Figure 10. Second Venus swing-by

**Application Example 2: Ariane case EVVMM strategy**

The Soyuz option allows to deliver about 600 kg of payload in orbit around Mercury with two swing-bys of Venus and two swing-bys of Mercury. The same strategy can be employed to deliver a higher mass but using a different launcher: an Ariane 5 with Vinci upper stage. Even in this case a 30 days coast arc is inserted before each swing-by and before the final encounter with Mercury.

Initial conditions are given by Ariane performances at launch while final conditions, as for the Soyuz case, depend on a final break with chemical propulsion.

Results for the Soyuz option and for the Ariane option are summarized in Table 3 where the date and the altitude for each swing-by is reported. All three cases have been analyzed taking into account Soyuz or Ariane 5 performances and a thrust dependent on the power provided by the solar arrays, apart from the first case for which the thrust has been kept constant for the transfer between the Earth and Venus and equal to half the maximum nominal thrust. A more realistic model leads to an improvement of the payload mass because the thrust increases, approaching Venus, and the first thrust arc reduces in length.

Table 3. Soyuz and Ariane 5 cases: summarizing table

| DATA | SOYUZ | | SOYUZ | | ARIANE 5 | |
|---|---|---|---|---|---|---|
| $v_\infty$ | 1.8 km/s | | 1.65 km/s | | 1.72 km/s | |
| Max Thrust | 0.17-0.34 | | 0.34 (power dependent) | | 0.6 (power dependent) | |
| Isp | 3200 s | | 3200 s | | 3200 s | |
| Launch Date | 23 July 2010 | | 21 July 2010 | | 22 July 2010 | |
| Initial Mass | 1316.8 kg | | 1334.4 kg | | 2629.0 kg | |
| Final Mass | 1030.1 kg | | 1042.3 kg | | 1984.9 kg | |
| Payload Mass | 589.4 kg | | 600.1 kg | | 1445.5 kg | |
| Arrival Velocity | 0.511 km/s | | 0.511 km/s | | 0.318 km/s | |
| Arrival Date | 27 Oct 2012 | | 27 Oct 2012 | | 10 Nov 2012 | |
| | Altitude | Date | Altitude | Date | Altitude | Date |
| Venus Swing-by | 3630 km | 28 Nov 2010 | 3629 km | 29 Nov 2010 | 5197 km | 02 Dec 2010 |
| Venus Swing-by | 300 km | 20 Mar 2011 | 300 km | 21 Mar 2011 | 300 km | 24 Mar 2011 |
| Mercury Swing-by | 200 km | 25 May 2012 | 200 km | 25 May 2012 | 200 km | 05 Jun 2012 |
| Mercury Swing-by | 200 km | 05 July 2012 | 200 km | 05 July 2012 | 200 km | 27 Jul 2012 |

**Application Example 3: Ariane case EVVMMM strategy**

An interesting alternative solution is to exploit resonance with Mercury to progressively reduce the perihelion. After the second swing-by of Venus the spacecraft is directly injected into a trajectory toward Mercury. The spacecraft encounters Mercury close to its perihelion and performs a first gravity maneuver acquiring a 2:3 resonant orbit synchronous with the motion of Mercury. After three rounds a second encounter reduces the aphelion leading to a 3:4 resonance with Mercury. Finally the spacecraft performs a third swing-by with Mercury and is injected into a trajectory so close to the orbit of Mercury that allows a relative arrival velocity of 340 m/s.

An Ariane5 at departure provides an asymptotic velocity of 2.8 km/s and the wet mass at launch is of 2320 kg. No special constraints are imposed on the declination but the solar arrays are dimensioned to allow a thrust of about 0.4 N up to Venus and then they operate at full power providing a thrust of 0.6 N down to Mercury. Not only has the thrust been taken variable but even the Isp has been considered equal to 3400 s at beginning of life and equal to 3200 s at end of life due to degradation of the efficiency of the SEP module. The obtained trajectory is represented in Figure 11, even in this case a circle marker represents departure from the Earth and a diamond represents arrival at Mercury. Stars represent swing-bys and, as expected, the first two swing-bys of Mercury occur close to the perihelion.

The effect of resonance with Mercury between swing bys is evident in Figures 12 and 13, reporting two significant orbital parameters. The semimajor axis benefits from the resonant condition, being progressively reduced, while orbit inclination is unaffected, as expected.

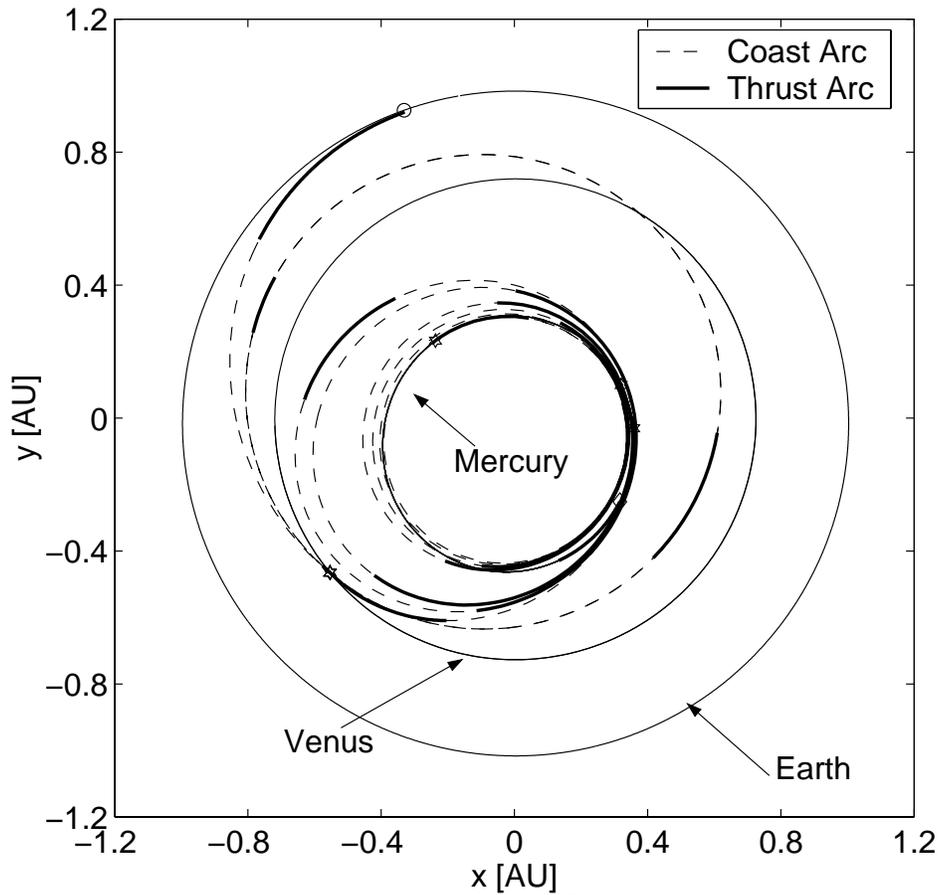

Figure 11. Trajectory in the ecliptic plane for the EVVMMM strategy

Table 4. Summarizing table for the EVVMMM strategy with Ariane 5

| DATA | ARIANE | |
|---|---|---|
| v∞ | 2.88 km/s | |
| Max Thrust | 0.4-0.6 | |
| Isp | 3400-3200 | |
| Launch Date | 09 Jan 2009 | |
| Initial Mass | 2320 kg | |
| Final Mass | 1882.1 kg | |
| Payload Mass | 1352.96 kg | |
| Arrival Velocity | 0.340 km/s | |
| Arrival Date | 01 Aug 2012 | |
| Venus Swing-by | 4217.1 km | 16 Apr 2009 |
| Venus Swing-by | 350 km | 10 July 2010 |
| Mercury Swing-by | 100 km | 6 Sep 2010 |
| Mercury Swing-by | 100 km | 1 Jun 2011 |
| Mercury Swing-by | 19280 km | 7 Jun 2012 |

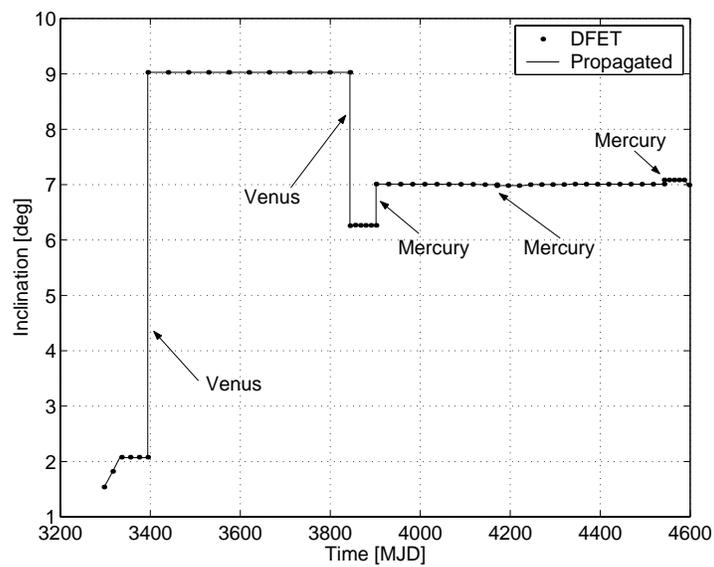

Figure 12. Time history of the semimajor axis

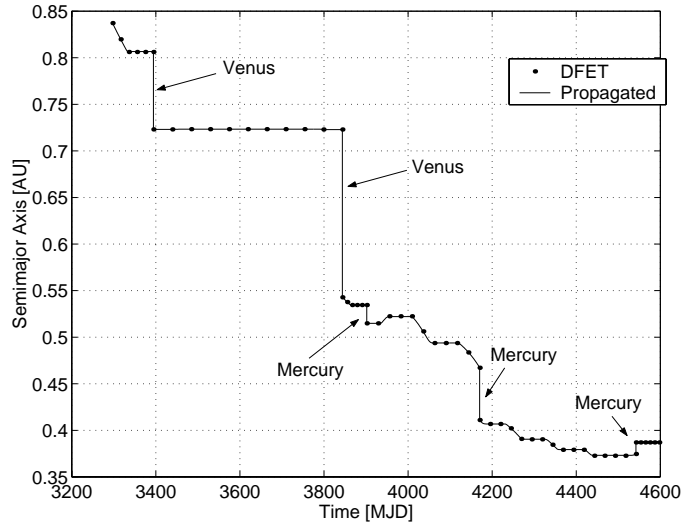

Figure 13. Time history of the inclination.

**CONCLUSIONS**

In this paper the problem of designing an optimal interplanetary transfer trajectory from the Earth to planet Mercury has been solved with a direct optimization approach and a transcription by Finite Elements in Time. The trajectory optimization problem is particularly complex due to the combination of low-thrust and multiple gravity assist maneuvers used to reduce the demands in terms of $\Delta v$. The problem is split into phases and for each one both states and controls are parameterized using DFET, an additional set of parameters is then included leading to a direct multiphase parametric optimization of the trajectory. Swing-bys are, at first, introduced through a simplified link-conic model for which the altitude is a parameter to be optimized then they are introduced as a full propagation of the hyperbolae. In the latter case orbital parameters of the hyperbolae are included among NLP parameters and optimized.

The parametric optimization using a combination of collocation by FET and shooting is quite robust and solves efficiently and accurately various problems with a reduced set of NLP variables. The effectiveness of the method is proved even by a comparison with an equivalent solution obtained with an indirect multiple shooting approach, on a sample mission to planet Mercury. For the mission optimization, two launch options and two different strategies are analyzed, in particular the EVVMMM strategy appears to be quite efficient in terms of $\Delta v$ gained from each swing-by.

**AKNOWLEDGEMENTS**

The authors would like to thank Dr. Rüdiger Jehn and Mr. Markus Katskowski of the European Space Operation Center (ESOC/ESA) for the data relative to the indirect solution and to Soyuz and Ariane performances.